\documentclass[a4paper]{article}

\usepackage{fullpage}

\usepackage{amssymb}
\usepackage{amsmath}
\usepackage{amsthm}
\usepackage{dsfont}

\usepackage{subfigure}

\usepackage{tikz}
\usetikzlibrary{arrows}
\usetikzlibrary{shapes.symbols,shapes.geometric,arrows,calc}
\usetikzlibrary{arrows,shapes,decorations,decorations.pathreplacing,decorations.pathmorphing,intersections,calc,fit}
\newcommand {\junk}[1]{}

\newtheorem{thm}{Theorem}
\newtheorem{lem}[thm]{Lemma}

\newtheorem{cor}[thm]{Corollary}

\newtheorem{ex}[thm]{Example}

\newcommand{\IR}{\mathbb{R}}
\newcommand{\IN}{\mathbb{N}}
\newcommand{\snorm}[1]{\delta\left({#1}\right)}
\newcommand{\graph}[1]{G\left({#1}\right)}
\newcommand{\onevec}{\mathbf{1}}

\def\ind{T}

\title{
Asymptotic Consensus Without Self-Confidence
}

\author{
Thomas Nowak\thanks{Thomas Nowak is with the D\'epartement d'Informatique, \'Ecole normale sup\'erieure, Paris, France. Email: {\tt thomas.nowak@ens.fr}}%
}

\date{}

\begin{document}

\maketitle
\thispagestyle{empty}
\pagestyle{empty}

\begin{abstract}
This paper studies asymptotic consensus in systems in which agents do not
necessarily have self-confidence, i.e., may disregard their own value during
execution of the update rule.
We show that the prevalent hypothesis of self-confidence in many convergence
results can be replaced by the existence of aperiodic cores.
These are stable aperiodic subgraphs, which allow to virtually store information
about an agent's value distributedly in the network.
Our results are applicable to systems with message delays and memory loss.
%
\end{abstract}



%
%
%
%
%
%
%
%
%

\section{Introduction}


Asymptotic consensus is a phenomenon observed in certain biological, physical,
and sociological systems.
It is also utilized in some engineered man-made computer systems.
The phenomenon consists in agents communicating in a very simple fashion to
asymptotically reach agreement on a common real value.
In nature, it can be observed (e.g.,\cite{Rey87,JLM03,TJP07,HK02})
in bird flocking, firefly synchronization,
synchronization of coupled oscillators, or opinion spreading.
In engineering, it is used for sensor fusion, dynamic load balancing protocols,
robot formation protocols, replication techniques, or rendezvous in space.

There is a very simple algorithm for asymptotic consensus that works in a
large class of environments:
In every computation step of a process, it updates its value to some average of
all values it has received, and then sends out its new value.
This simple algorithm has two remarkable properties:
Firstly, it is very simple and yet manages to solve asymptotic consensus in a
surprisingly large number of different environments.
Secondly, it is an algorithm that can be {\em observed\/} in nature.
More specifically, it serves as a widely accepted model in biology, physics,
and sociology to explain various phenomena such as bird flocking,
synchronization of coupled oscillators, and opinion spreading.
It thus stands to reason to expect the algorithm to have a certain robustness
against adverse environments.
Of course, one can think of using it to attain approximate agreement in
man-made, engineered, systems.
And indeed, it is actually used, for example, in sensor fusion.
For engineered systems, the viewpoint is not one of observing and explaining a
given system, but of {\em analyzing\/} it for prediction of its future behavior
or for assessing the need to improve the system.
The speed of convergence in the context of asymptotic consensus is a measure for
the stabilization time, or the transient phase, of the system.
Obviously, the sharper the analysis of the system and its performance, the
tighter it can be integrated into the timing constraints of a larger system,
and hence the larger the potential performance of the larger system.

The analysis becomes significantly harder if the communication graphs, or the
weights, change over time, if communication delays are introduced and if nodes
are susceptible to certain faults.
If one admits the dynamicity of the communication graph, then one has already
accounted for a large class of faults, namely link faults.
The addition of communication delays covers timing faults on links.
A class of faults that has received considerably less attention in the
literature is that of memory faults, either by memory loss or memory delays,
i.e., the value read from local memory is not that of the most recent write
operation.
Memory delays become more probable with the advent of modern pipelined
architectures and memories with weakened consistency properties.
The present paper has as its goal the study of systems in which processes
cannot, or do not, access their most recent value, but may read an older one or
disregard it altogether.
In the context of natural asymptotic consensus systems like in sociology, this
phenomenon is more naturally called a {\em lack of self-confidence\/} and has
its specific interest in the analysis of such systems.
The paper extends a variety of convergence results known for cases with
self-confidence to cases without and identifies the importance of having a
certain replacement for self-confidence, which we call aperiodic cores.
Self-confidence is a specific instance of this notion.
Moreover, we discuss an explicit example showing the boundary between
convergence and non-convergence in the context of aperiodic cores, shedding a more precise light on the frontier.

In linear algebraic terms, the study of asymptotic consensus is the study of
infinite backwards products of stochastic matrices.
The first convergence result for products of stochastic matrices is the
Perron-Frobenius theorem, which states that the powers of an
ergodic stochastic matrix converge to a rank~$1$ stochastic matrix.
%
It was first generalized to a non-constant product of matrices by
Wolfowitz~\cite{Wol63} who showed that if every {\em finite\/} product of matrices of a set~$\mathcal{M}$ of matrices is ergodic, then
every backwards of matrices in~$\mathcal{M}$ converges to a rank~$1$
stochastic matrix.
The strict finiteness and ergodicity conditions in Wolfowitz' theorem were found
to be inappropriate for many applications.
Subsequently, Wolfowitz' theorem was extended in several directions
(see, for example, \cite{BHOT05}, \cite[Section II.G]{OSFM07},
or \cite{CMA08a}).
However, no direct generalization of Wolfowitz' theorem or the Perron-Frobenius
theorem was obtained.
This is due to the fact that these results all assume a strictly positive
diagonal in all occurring matrices.
In this sense, the results on asymptotic consensus in dynamic settings are no strict generalizations of the Perron-Frobenius theorem or Wolfowitz' theorem, precisely because of the fact that they require a strictly  positive diagonal.
One goal of this paper is to remedy this deficiency; by providing convergence results for asymptotic consensus in dynamic settings without this hypothesis.
Thus, our results are both strict generalizations of the Perron-Frobenius
theorem and existing convergence theorems in asymptotic consensus.

The rest of the paper is organized as follows:
Section~\ref{sec:prelim} introduces the model, discusses related work, and
gives necessary preliminary results.
The notion of aperiodic cores is defined in Section~\ref{sec:cores} and the
first new convergence result based on this notion follows in
Section~\ref{sec:coordinated}.
We generalize the definition of aperiodic cores in
Section~\ref{sec:clusterings} by introducing the notion of clusterings.
This is useful to talk about hierarchic systems with local leader agents, as
they naturally appear in the reduction from non-synchronous to synchronous
settings.
We apply this notion in Sections~\ref{sec:dyn:coord}, \ref{sec:dyn:fixed},
and~\ref{sec:red} to show quite general convergence theorems in various
environments, together with upper bounds on the convergence rate where
applicable.
Each of our theorems is followed by a corollary in form of an already known
result in the literature.
We do this to facilitate finding the context in terms of classical results in
which the present paper generalizes the state of the art.
Section~\ref{sec:conclusion} concludes the paper with some final remarks.

\section{Asymptotic Consensus}\label{sec:prelim}

In particular in computer science for multi-agent
systems whose agents start with a private value and repeatedly form averages of
perceived values of others.
These types of multi-agent systems are not only used in computer networks, but
have also been found to model various physical and
biological phenomena like the behavior of bird flocks~\cite{BHOT05,OSFM07}.
Mathematically, they translate into long and infinite backwards products of
stochastic matrices.

\subsection{Computational Model}

The distributed computing model in which we study asymptotic consensus is the
following:
There are~$n$ distinguishable agents, each agent $i\in[n] = \{1,2,\dots,n\}$
possessing a real state
variable~$x_i$ and communicating by exchanging messages.
There is a global discrete time base, referred to by nonnegative integers
in $\IN=\{0,1,2,\dots\}$.
At every time $t\in\IN$, we denote the content of the agents' state variables
by $x_i(t)$.
The initial value of state variable~$x_i$ is $x_i(0)$.
At every time $t\in\IN$, every agent sends the content of its state variable to
all other agents.
Messages may be delayed and/or lost.
All agents simultaneously update their state variable at all positive times
$t=1,2,3,\dots$ 
to some weighted average value of the received values, at most one of each
other agent,
and its current content of its own state variable.

Since the new content of the state variable is a mean value, 
there
exists a \mbox{$\Delta_{i,j}(t)>0$} for every $j\in[n]$ such that
\begin{equation}\label{eq:evolution:delays}
x_i(t) = \sum_{j=1}^n A_{i,j}(t) \cdot x_j\big(t-\Delta_{i,j}(t)\big)
\end{equation}
with
\begin{equation}\label{eq:evolution:delays:average}
\sum_{j=1}^n A_{i,j}(t) = 1
\enspace.
\end{equation}
A {\em configuration\/} of asymptotic consensus is a collection of real values,
one for each agent's state variable, i.e., a vector in~$\IR^n$.
An {\em execution\/} of asymptotic consensus is an infinite sequence of
configurations $x(t)\in\IR^n$ following the
evolution~\eqref{eq:evolution:delays} for some choice of the~$A_{i,j}(t)$
and the~$\Delta_{i,j}(t)$.
An execution {\em reaches asymptotic consensus\/} if~$x(t)$ converges and all
component-wise limits $\lim_{t\to\infty}x_i(t)$ are equal.

We call an {\em averaging matrix\/} a matrix whose entries are all nonnegative
and whose row sums are all~$1$.
In other words, it is a row stochastic matrix.
Equation~\eqref{eq:evolution:delays:average} assures that the collection of
the~$A_{i,j}(t)$ is an averaging matrix for all~$t$.
A {\em delay matrix\/} for time~$t$ is a matrix of integers between~$1$
and~$t$.
For every~$t$, the collection of the~$\Delta_{i,j}(t)$ is a delay matrix
for~$t$.
Hence an execution is determined by the initial configuration~$x(0)$,
the sequence of the averaging matrices~$A(t)$, and the sequence of the delay
matrices~$\Delta(t)$.
A pair consisting of a sequence of averaging
matrices~$A(t)$ and a
sequence of vectors~$\Delta(t)$ such that every~$\Delta(t)$ is a delay matrix
for~$t$ is referred to as a {\em setting}. 
An {\em environment\/} is a nonempty set of settings.
We say that a setting or an environment reaches asymptotic consensus if all of
its executions do.

An important parameter of a setting is the maximum entry of the delay matrices,
if it exists.
We call a setting {\em $B$-bounded\/} if all entries of its delay matrices are
at most~$B$.
A $1$-bounded setting is called {\em synchronous\/} and is determined uniquely
by the sequence of averaging matrices.
If the nonzero entries of the averaging matrices are lower bounded by some
positive~$\alpha$, then we say that the setting has {\em minimal confidence~$\alpha$}.
It has {\em self-confidence\/} if all diagonal entries are positive.
The {\em communication digraph\/} of a stochastic matrix~$A$ in $\IR^{n\times
n}$ has node set~$[n]$ and contains an edge $(i,j)$ if and only if $A_{i,j}>0$.

We note that not every non-synchronous setting reaches asymptotic consensus;
{\em not even\/} with self-confidence and strongly connected bidirectional
communication graphs.
The following example shows this.
The problem arises if the delay $\Delta_{i,i}(t)$ is strictly greater than~$1$,
i.e., node~$i$ does not use its most recent value for the update rule.
It is one of the goals of the present paper to study sufficient conditions that
enable convergence even if $\Delta_{i,i}(t)>1$ for some, or even all,~$i$
and~$t$.

\begin{ex}\label{ex:period}
With $n=2$ agents, we choose the averaging matrices
\begin{equation*}
A(1) = 
\begin{pmatrix}
1 & 0\\
0 & 1
\end{pmatrix}
\quad
\text{and}
\quad
A(t) = 
\begin{pmatrix}
1/2 & 1/2\\
1/2 & 1/2
\end{pmatrix}
\text{for }
t\geq 2
\end{equation*}
and the initial vector $x(0) = {}^t(0, 1)$.
Thus there is self-confidence and a minimal confidence of $1/2$.
For the delay matrices, we choose
\begin{equation*}
\Delta(1) = 
\begin{pmatrix}
1 & 1\\
1 & 1
\end{pmatrix}
\quad
\text{and}
\quad
\Delta(t) = 
\begin{pmatrix}
2 & 1\\
1 & 2
\end{pmatrix}
\text{for }
t\geq 2
\enspace,
\end{equation*}
i.e., for times $t\geq 2$, there is a delay to itself at every agent of~$2$
(even though the delay to the other agent is~$1$).
The communication graph for $t\geq 2$ is shown in Fig.~\ref{fig:ex:period:1}.
One can show that $x_1(2t)\to 1/3$ as $t\to\infty$ whereas $x_1(2t+1)\to 2/3$.
Similarly, $x_2(2t)\to 2/3$ and $x_2(2t+1) \to 1/3$.
That is, the system is asymptotically periodic with period~$2$.
The issue becomes clearer when looking at the equivalent synchronous system as studied by Cao, Morse, and Anderson~\cite{CMA08b}.
Its communication graph for $t\geq2$ is depicted in Fig~\ref{fig:ex:period:2}.
This equivalent synchronous communication graph has a period of~$2$.
We introduce their reduction in more detail at the end of
Section~\ref{sec:related}.
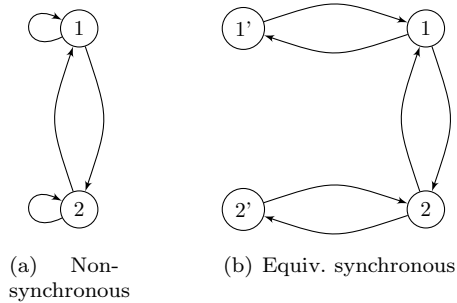
\begin{figure}
\centering
\subfigure[Non-synchronous]
{
\label{fig:ex:period:1}
\begin{tikzpicture}[>=latex',scale=0.8, every node/.style={transform shape}]
\node[draw,circle] (n2) at (3,0) {1};
\node[draw,circle] (n3) at (3,-3) {2};
\draw[->] (n2) .. controls +(0.5,-1.5) .. (n3);
\draw[->] (n3) .. controls +(-0.5,1.5) .. (n2);
\draw[->] (n3) .. controls +(-1.0,-0.5) and +(-1.0,0.5) .. (n3);
\draw[->] (n2) .. controls +(-1.0,-0.5) and +(-1.0,0.5) .. (n2);
\end{tikzpicture}
}
\hspace*{1cm}
\subfigure[Equiv.\ synchronous]
{
\label{fig:ex:period:2}
\begin{tikzpicture}[>=latex',scale=0.8, every node/.style={transform shape}]
\node[draw,circle] (n1) at (0,0) {1'};
\node[draw,circle] (n2) at (3,0) {1};
\node[draw,circle] (n3) at (3,-3) {2};
\node[draw,circle] (n4) at (0,-3) {2'};
\draw[->] (n1) .. controls +(1.5,0.5) .. (n2);
\draw[->] (n2) .. controls +(-1.5,-0.5) .. (n1);
\draw[->] (n2) .. controls +(0.5,-1.5) .. (n3);
\draw[->] (n3) .. controls +(-0.5,1.5) .. (n2);
\draw[->] (n4) .. controls +(1.5,0.5) .. (n3);
\draw[->] (n3) .. controls +(-1.5,-0.5) .. (n4);
\end{tikzpicture}
}
\caption{Communication graphs for $t\geq 2$ in the original non-synchronous and the equivalent synchronous setting in Example~\ref{ex:period}}
\end{figure}
\end{ex}

In a synchronous setting, the evolution of configurations~$x(t)$ is
governed by the linear recursive law
\begin{equation*}
x(t) = A(t) \cdot x(t-1)
\end{equation*}
where~$A(t)$ is a row stochastic matrix.
Defining the product matrices
\begin{equation*}
P(t) = A(t)\cdot A(t-1) \cdots A(1)
\enspace,
\end{equation*}
we have 
\begin{equation*}
x(t) = P(t) \cdot x(0)
\enspace. 
\end{equation*}
In particular, the sequence of state vectors is
determined by the initial vector~$x(0)$ and the sequence of row stochastic
matrices~$A(t)$.

In the following sections, we will also use the notation
\begin{equation*}
P(t,s) = A(t) \cdot A(t-1) \cdots A(s+1)
\end{equation*}
for partial products.
It is $P(t) = P(t,0)$ for all~$t$ and $P(t,s) = I$, the identity matrix, if
$t\leq s$.
If all $A(t)$ are equal to a constant matrix~$A$, then $P(t) = A^t$.

%
%

\subsection{Related Work}\label{sec:related}
In this subsection, we list several convergence theorems in the literature that
our results generalize.
All of them suppose self-confidence.

Tsitsiklis introduced the {\em bounded intercommunication\/} assumption.
It states that if an edge $(i,j)$ appears in infinitely many communication
digraphs, then is appears in one of the digraphs $$\graph{A(t)},
\graph{A(t+1)}, \dots, \graph{A(t+B-1)}$$ for a fixed~$B$ and all~$t$.

\begin{thm}[Tsitsiklis \cite{Tsi84}]
\label{thm:tsitsiklis}
A synchronous setting with averaging matrices $A(1),A(2),\dots$
with self-confidence and minimal confidence~$\alpha$ reaches asymptotic
consensus if the digraph~$G_\infty$ formed by the edges appearing in infinitely
many communication digraphs is strongly connected and the bounded
intercommunication assumption holds.
\end{thm}

Moreau and Hendrickx and Blondel independently showed that the bounded
intercommunication  assumption can be replaced by the assumption that every
communication digraph is bi-directional:

\begin{thm}[Moreau \cite{Mor05}, Hendrickx and Blondel \cite{HB05}]
\label{thm:moreau}
A synchronous setting with averaging matrices $A(1),A(2),\dots$
with self-confidence and minimal confidence~$\alpha$ reaches asymptotic
consensus if the digraph~$G_\infty$ 
is strongly connected and every communication
digraph is bi-directional.
\end{thm}

Blondel et al.\ generalized this result to $B$-bounded settings:

\begin{thm}[Blondel et al.~\cite{BHOT05}]
\label{thm:blondel:et:al}
A $B$-bounded setting with averaging matrices $A(1),A(2),\dots$
with self-confidence and minimal confidence~$\alpha$ reaches asymptotic
consensus if the digraph~$G_\infty$ 
 is strongly connected and every communication
digraph is bi-directional.
\end{thm}

Touri and Nedi\'c generalized the assumption of bi-directional
digraphs to
digraphs that are completely reducible.
Charron-Bost recently showed its extension to $B$-bounded
settings.

\begin{thm}[Touri and Nedi\'c~\cite{TN11}, Charron-Bost~\cite{CB13}]
\label{thm:touri:nedic}
A $B$-bounded setting with averaging matrices $A(t)$
with self-confidence and minimal confidence~$\alpha$ reaches asymptotic
consensus if the digraph~$G_\infty$ 
is strongly connected and every communication
digraph is completely reducible.
\end{thm}

If an execution $x(t)$ reaches asymptotic consensus, one can ask the question
of the speed at which this convergence occurs.
Olshevsky and Tsitsiklis noted that this speed tends to be exponential and
have hence defined the {\em rate of convergence\/} as
\begin{equation*}
\lim_{t\to\infty}\lVert x(t) - x^* \rVert_2^{1/t}
\enspace.
\end{equation*}

Cao, Morse, and Anderson studied {\em coordinated\/} communication digraphs,
i.e., digraphs that have a node~$j$ such that every other node has a path
to~$j$.
They obtained the following result:

\begin{thm}[Cao, Morse, and Anderson \cite{CMA08a,CMA08b}]
\label{thm:cao:et:al}
A $B$-bounded setting with sequence of  averaging matrices $A(1),A(2),\dots$
with self-confidence and minimal confidence~$\alpha$ reaches asymptotic
if every communication digraph is coordinated.
Moreover, the rate of convergence is less than~$1$.
\end{thm}

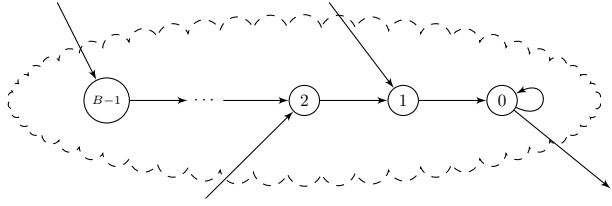
\begin{figure}
\centering
\begin{tikzpicture}[>=latex',scale=0.65, every node/.style={transform shape}]
\node[draw,circle] (n0) at (0,0) {$0$};
\node[draw,circle] (n1) at (-2,0) {$1$};
\node[draw,circle] (n2) at (-4,0) {$2$};
\node (n3) at (-6,0) {$\cdots$};
\node[draw,circle] (n4) at (-8,0) {$\scriptstyle B-1$};
\draw[->] (n1) -- (n0);
\draw[->] (n2) -- (n1);
\draw[->] (n3) -- (n2);
\draw[->] (n4) -- (n3);
\draw[->] (-9,2) -- (n4);
\draw[->] (-6,-2) -- (n2);
\draw[->] (-3.5,2) -- (n1);
\draw[->] (n0) .. controls +(1.0,-0.5) and +(1.0,0.5) .. (n0);
\draw[->] (n0) -- (2.2,-1.8);
\node[dashed,cloud, cloud puffs=54, draw,minimum width=12cm, minimum
height=3.5cm] at (-4,0) {};
\end{tikzpicture}
\caption{The $B$ copies of an agent in Cao, Morse, and Anderson's reduction}
\label{fig:delay:reduction}
\end{figure}

To prove their result, they described a reduction of $B$-bounded settings to
synchronous
settings, albeit with~$B$ times as many agents as the original
setting~\cite[Section 4.1]{CMA08b}.
The idea is to replicate every agent~$B$ times, but to shift the copies in
time, i.e., at time~$t$ there is one copy holding the value~$x_i(t)$, one
$x_i(t-1)$, and so on until $x_i(t-B+1)$.
This results in synchronous setting for asymptotic consensus.
The replication of agents is illustrated in Fig.~\ref{fig:delay:reduction}.
Only the copy for the current value $x_i(t)$ has links to other agents' copies.
Nonetheless, no such restriction exists for incoming edges.
In the new resulting communication digraphs, even if all agents have self-loops
in the original communication digraphs, not all nodes have them.

\subsection{Dobrushin Semi-Norm for Stochastic Matrices}

All stochastic matrices have~$1$ as an eigenvalue of maximum modulus.
If the matrix is irreducible, the corresponding right-eigenspace is
one-dimensional
and generated by the column vector $\mathbf{1} = {}^t(1,1,\dots,1)$.
When studying such matrices, we are hence led to consider the distance of
vectors to this eigenspace.
Indeed, we will see that considering this
distance is an appropriate tool for products of stochastic matrices.



The Dobrushin vector semi-norm on~$\mathds{R}^n$ is defined by setting
\[
\delta( x) = \inf_{y\in
\mathds{R}\cdot \mathbf{1}} \lVert x-y \rVert_\infty
\enspace.
\]
%
This vector semi-norm induces the Dobrushin matrix semi-norm
on~$\mathds{R}^{n\times n}$ by
defining it in the operator norm fashion:
\[
\delta( A ) = \sup_{\substack{x\in\mathds{R}^n\\ \delta( x) \neq 0}}
\frac{\delta( Ax )}{\delta(x)}
\]
Clearly, $\delta( A ) = 0$ if the image of~$A$ is contained in the
subspace~$\mathds{R}\cdot\mathbf{1}$.

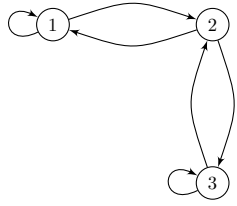
\begin{figure}
\centering
\begin{tikzpicture}[>=latex',scale=0.7, every node/.style={transform shape}]
\node[draw,circle] (n1) at (0,0) {1};
\node[draw,circle] (n2) at (3,0) {2};
\node[draw,circle] (n3) at (3,-3) {3};
\draw[->] (n1) .. controls +(1.5,0.5) .. (n2);
\draw[->] (n2) .. controls +(-1.5,-0.5) .. (n1);
\draw[->] (n2) .. controls +(0.5,-1.5) .. (n3);
\draw[->] (n3) .. controls +(-0.5,1.5) .. (n2);
\draw[->] (n1) .. controls +(-1.0,-0.5) and +(-1.0,0.5) .. (n1);
\draw[->] (n3) .. controls +(-1.0,-0.5) and +(-1.0,0.5) .. (n3);
\end{tikzpicture}
\caption{Digraph~$\graph{A}$ of matrix~$A$}
\label{fig:ex:snorm:alpha}
\end{figure}

We now give an example of a matrix whose semi-norm is strictly less than~$1$,
but that has neither a strictly positive column
nor a strictly positive diagonal.
The matrix is equal to
\begin{equation*}
A = 
\begin{pmatrix}
1/2 & 1/2 & 0 \\
1/2 & 0 & 1/2 \\
0 & 1/2 & 1/2
\end{pmatrix}
\end{equation*}
and its digraph is depicted in Fig.~\ref{fig:ex:snorm:alpha}.
In fact, $\snorm{A}$ is equal to~$1/2$.


The following lemma characterizes the matrices with a Dobrushin semi-norm
strictly smaller than~$1$.
It uses the notion of a {\em scrambling\/} matrix.
A stochastic matrix~$A$ is scrambling if for all indices~$i_1,i_2$ 
there exists an index~$j$ such that both $A_{i_1,j} >0$ and $A_{i_2,j} > 0$.
Note that, a fortiori, $A$ is scrambling if it has a strictly positive column.
Its proof follows from the formula 
\begin{equation*}
\delta({A}) = \max_{i_1,i_2\in[n]}\sum_{j=1}^n (A_{i_1,j} -
A_{i_2,j})_+
\end{equation*}
for the Dobrushin matrix semi-norm
where $(x)_+ = \max\{x,0\}$ denotes the positive part of~$x$.

\begin{lem}[\cite{CMA08a,CB13}]\label{lem:norm:bound}
Let~$A$ be a stochastic matrix.
We always have $\delta( A ) \leq 1$
and  $\delta(A) < 1$ if and only if~$A$ is scrambling.
In this case, $\delta(A) \leq 1 - \alpha$ where~$\alpha$ is the smallest
nonzero entry of~$A$.
\end{lem}
%

The next lemma shows the utility of~$\delta$ to show convergence and asymptotic
agreement.

\begin{lem}\label{lem:norm:limit}
The sequence of partial backwards products~$P(t)$ converges to a
rank~$1$ stochastic matrix if and only if $\delta\big( P(t)\big)\to 0$ as
$t\to\infty$.
\end{lem}
\begin{proof}
If $P(t)$ converges to a rank~$1$ stochastic matrix~$P$, then $\delta(P)=0$.
By continuity of~$\delta$ and monotonicity of~$\delta\big(P(t)\big)$,
necessarily $\delta\big(P(t)\big) \to 0$.

To prove the converse implication, we now assume that $\delta\big(P(t)\big) \to
0$.
We show that, for every $x\in\mathds{R}^n$, the sequence of vectors $P(t)\cdot
x$ converges by showing that it is Cauchy.
This then concludes the proof because stochasticity is preserved when taking
limits and~$\delta$ is continuous.

Let $\varepsilon>0$.
Because also $\delta\big( P(t)\cdot x \big) \to 0$, there exists some~$T$ such
that  $\delta\big( P(T) \cdot x\big) \leqslant
\varepsilon/2$.
Letting $y\in\mathds{R}\cdot \mathbf{1}$ such that $\delta\big( P(T)\cdot
x\big) =
\lVert P(T)\cdot x - y \rVert_\infty$, we calculate for every $t\geqslant T$:
\begin{align}
\lVert P(t) \cdot x - P(T) \cdot x\rVert_\infty & \leqslant 
\lVert P(t,T)\cdot P(T) \cdot x - y \rVert_\infty\notag\\
&\quad + \lVert P(T) \cdot x -
y\rVert_\infty
\notag\\ & =
\lVert P(t,T)\cdot ( P(T) \cdot x - y) \rVert_\infty\notag\\
&\quad + \lVert P(T) \cdot x -
y\rVert_\infty
\notag
\\ & \leqslant 
2 \cdot \lVert P(T) \cdot x -
y\rVert_\infty
\notag\\
& =
2\cdot \delta\big( P(T)\cdot x\big)
\leqslant \varepsilon
\notag
\end{align}
because $y = P(t,T)\cdot y$ since $P(t,T)$ is stochastic and~$y$ is a multiple
of~$\onevec$.
This shows that $P(t)\cdot x$ is indeed a Cauchy sequence.
\end{proof}

We now provide a tool to prove convergence of the matrix semi-norm of
a product to zero by stating a sufficient condition for the semi-norm of a
factor to be constantly bounded away from~$1$.
It shows in particular that the semi-norm of a stochastic matrix is at
most~$1$.

\subsection{Graph Interpretation of Matrix Products}

Let~$i$ and~$j$ be nodes of a digraph~$G$.
A {\em walk\/} in~$G$ from~$i$ to~$j$ is a finite sequence of adjacent nodes
in~$G$ that starts at~$i$ and ends at~$j$.
Its length is the number of nodes in the sequence minus one.

The following lemma characterizes positivity of entries in products of
stochastic matrices solely in terms of the matrices' associated digraphs.
It should be noted that, because we study backward products, the walks grow at
the start node and not at the end node.

\begin{lem}\label{lem:graph}
Let $0\leq s\leq t$ and $i,j\in[n]$.
Then $P_{i,j}(t,s)$ is positive if and only if there exist
$i_t,i_{t-1},\dots,i_s\in[n]$ with $i_t=i$ and $i_s=j$ such
that $(i_\tau , i_{\tau-1})$ is an
edge of~$G\big(A(\tau)\big)$ for all $s+1\leq  \tau\leq t$.
\end{lem}

If a strongly connected digraph is aperiodic, there exist walks of arbitrary
length between all pairs of nodes as long as the length is greater or equal to
a number called the {\em exponent\/} (sometimes also {\em index}\/) of the
digraph.
Formally, we denote the smallest~$T$ such that there is a walk from~$i$ to~$j$
of length~$t$ for all nodes~$i$ and~$j$ such that~$j$ is reachable from~$i$
in~$G$  and all $t\geq T$ by~$T(G)$.
Wielandt provided an upper bound on the exponent, although many more
followed~\cite{DM64, Sch70, Kim79, GKP95, MNSS14}.
Wielandt's bound is the best possible upper bound in terms of only the number
of nodes.
If other parameters of the graph are known, however, tighter bounds exist.
Since the exponent~$T(G)$ appears in some of our bounds, it may be worthwhile 
to find a more precise bound for the specific graph appearing in a given
application framework.

\begin{thm}[Wielandt~\cite{Wie50}]\label{thm:wielandt}
Let~$G$ be a strongly connected aperiodic digraph with~$n$ nodes
Then the exponent of~$G$ is bounded by
\begin{equation*}
T(G) \leq W(n) = \begin{cases}
n^2 - 2n + 2 & \text{if } n\geq 2\\
0 & \text{if } n=1
\enspace.
\end{cases}
\end{equation*}
\end{thm}

\section{Aperiodic Cores}\label{sec:cores}

Classically, in asymptotic consensus, self-confidence of the agents is assumed.
That is, every communication digraph contains self-loops at all nodes.
This can model the fact that an agent does not ignore or forget its own
previous value.
We generalize the existence of self-loops, however:
A missing self-loop in a specific communication digraph can model memory loss
of an agent.
We replace the assumption of self-loops to {\em aperiodic cores}, which are
sub-digraphs of all of the settings' communication digraphs.
They can be seen as a ``distributed safety net against memory loss''.
In this sense, existence of self-loops is the assumption of a non-distributed
safety measure against memory loss or temporary self-distrust.
Their function in the proofs is similar to that of self-loops, but they are
more general.
A parameter that we use over and over in our results is that of the 
exponent of the aperiodic core.
If one assumes self-loops, then~$H$ only consists of self-loops at all nodes
and this parameter is equal to~$0$.
So, in our theorem statements, if one assumes self-confidence, then
$\ind(H)=0$.

We call a node~$j$ in a digraph~$G$ a {\em leader\/} of another node~$i$ if~$G$
contains a path from~$i$ to~$j$.
A digraph is {\em $j$-coordinated\/} if~$j$ is a leader of every node.
In this case, node~$j$ is called a {\em leader\/} of~$G$.
A digraph is {\em coordinated\/} if it is $j$-coordinated for some~$j$.
If~$j$ is a node of a digraph~$G$, we say that~$G$ is {\em $j$-aperiodic\/} if
$j$'s strongly connected component in~$G$ is primitive. 
A digraph~$H$ is a {\em core\/} of a sequence $G_1,G_2,\dots$ of
digraphs if~$H$ is a sub-digraph of every~$G_t$.

%
%

\section{Coordinated Aperiodic Cores}\label{sec:coordinated}

We start with assuming that there is a core that is coordinated and
leader-aperiodic.
The assumption of a core in particular applies if the communication digraph is
constant.
We hence get a direct generalization of the constant ergodic case:

\begin{thm}\label{thm:coord:const}
A synchronous setting with averaging matrices $A(t)$ 
with spanning core~$H$ and minimal
confidence~$\alpha$ reaches asymptotic consensus if
there exists some agent~$j_0$ such that~$H$ is $j_0$-coordinated and
$j_0$-aperiodic.
Moreover, the rate of convergence is at most
$1 - \alpha^{\ind(H)}/\ind(H)$.
\end{thm}

%
%

\begin{figure}
\centering
\begin{tikzpicture}[>=latex',scale=0.5, every node/.style={transform shape}]
\node[draw,circle] (n1) at (0,0) {1'};
\node[draw,circle] (n2) at (3,0) {1};
\node[draw,circle] (n3) at (3,-3) {2};
\node[draw,circle] (n4) at (0,-3) {2'};
\draw[->] (n1) .. controls +(1.5,0.5) .. (n2);
\draw[->] (n2) .. controls +(-1.5,-0.5) .. (n1);
\draw[->] (n2) -- (n3);
\draw[->] (n3) -- (n1);
\draw[->] (n4) .. controls +(1.5,0.5) .. (n3);
\draw[->] (n3) .. controls +(-1.5,-0.5) .. (n4);
\end{tikzpicture}
\caption{Variant of Example~\ref{ex:period} that converges}
\label{fig:ex:period:3}
\end{figure}
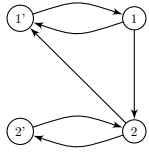

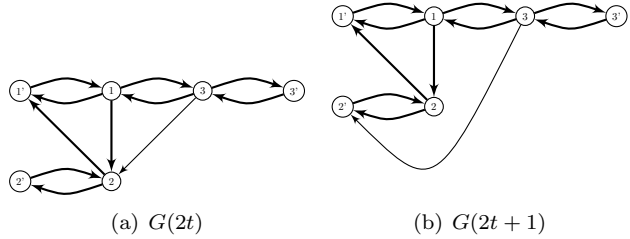
\begin{figure}
\centering
\subfigure[$G(2t)$]{
\begin{tikzpicture}[>=latex',scale=0.4, every node/.style={transform shape}]
\node[draw,circle] (n1) at (0,0) {1'};
\node[draw,circle] (n2) at (3,0) {1};
\node[draw,circle] (n3) at (3,-3) {2};
\node[draw,circle] (n4) at (0,-3) {2'};
\node[draw,circle] (n5) at (6,0) {3};
\node[draw,circle] (n6) at (9,0) {3'};
\draw[->,thick] (n1) .. controls +(1.5,0.5) .. (n2);
\draw[->,thick] (n2) .. controls +(-1.5,-0.5) .. (n1);
\draw[->,thick] (n2) -- (n3);
\draw[->,thick] (n3) -- (n1);
\draw[->,thick] (n4) .. controls +(1.5,0.5) .. (n3);
\draw[->,thick] (n3) .. controls +(-1.5,-0.5) .. (n4);
\draw[->,thick] (n2) .. controls +(1.5,0.5) .. (n5);
\draw[->,thick] (n5) .. controls +(-1.5,-0.5) .. (n2);
\draw[->,thick] (n5) .. controls +(1.5,0.5) .. (n6);
\draw[->,thick] (n6) .. controls +(-1.5,-0.5) .. (n5);
\draw[->] (n5) -- (n3);
\end{tikzpicture}
}
\subfigure[$G(2t+1)$]{
\begin{tikzpicture}[>=latex',scale=0.4, every node/.style={transform shape}]
\node[draw,circle] (n1) at (0,0) {1'};
\node[draw,circle] (n2) at (3,0) {1};
\node[draw,circle] (n3) at (3,-3) {2};
\node[draw,circle] (n4) at (0,-3) {2'};
\node[draw,circle] (n5) at (6,0) {3};
\node[draw,circle] (n6) at (9,0) {3'};
\draw[->,thick] (n1) .. controls +(1.5,0.5) .. (n2);
\draw[->,thick] (n2) .. controls +(-1.5,-0.5) .. (n1);
\draw[->,thick] (n2) -- (n3);
\draw[->,thick] (n3) -- (n1);
\draw[->,thick] (n4) .. controls +(1.5,0.5) .. (n3);
\draw[->,thick] (n3) .. controls +(-1.5,-0.5) .. (n4);
\draw[->,thick] (n2) .. controls +(1.5,0.5) .. (n5);
\draw[->,thick] (n5) .. controls +(-1.5,-0.5) .. (n2);
\draw[->,thick] (n5) .. controls +(1.5,0.5) .. (n6);
\draw[->,thick] (n6) .. controls +(-1.5,-0.5) .. (n5);
\draw[->] (n5) .. controls +(-3,-6) .. (n4);
\end{tikzpicture}
}
\caption{Equivalent synchronous communication graphs that alternate in time}
\label{fig:ex:period:4}
\end{figure}

We remark that Theorem~\ref{thm:coord:const} in particular shows that the
setting of Example~\ref{ex:period} reaches asymptotic consensus if we
change the delay $\Delta_{2,1}(t) = 2$, i.e., {\em increase\/} the message delay
from agent~$1$ to agent~$2$, for $t\geq 2$.
Indeed, the
resulting equivalent synchronous setting has an aperiodic core from time $t=2$
on, as is shown in Fig.~\ref{fig:ex:period:3}.
Note that, as the resulting stochastic matrix for the synchronous system is
ergodic and constant, that also the Perron-Frobenius theorem shows convergence
to asymptotic consensus in this case.
However, embedding this structure into a slightly larger but simple system
of~$3$ agents, as in
Fig.~\ref{fig:ex:period:4} (the aperiodic core is almost the whole graph and is shown in
bold; only a single edge changes continuously over time) shows the need the
generalization that Theorem~\ref{thm:coord:const} provides.

We prove this theorem in the rest of the subsection.

In general,
given a sequence of stochastic matrices $A(1),A(2),\dots$ in $\IR^{n\times n}$
and a node~$j\in[n]$, we 
define $S_{j}(t,s)$ to be the set of indices~$i\in[n]$ such
that~$P_{i,j}(t,s)$ is positive.
Denote by~$\mu_{j}(t,s)$ the smallest (positive) $P_{i,j}(t,s)$ with $i\in
S_{j}(t,s)$.
We also define $S_j(t) = S_j(t,0)$ and $\mu_j(t) = \mu_j(t,0)$.

It is easy to see that
$\mu_{j}(t,s) \geq \alpha^{t-s}$
if~$\alpha$ is the minimal confidence.
This will be our main tool to bound the convergence rate:
If $S_j(t,s) = [n]$, then $\snorm{P(t,s)} \leq 1 - \alpha^{t-s}$ by
Lemma~\ref{lem:norm:bound}.
And if we can show $S_j(t,s) = [n]$ whenever $t-s\geq T$ where~$T$ is some
constant, then 
\begin{equation*}
\begin{split}
\lim_{t\to\infty} \snorm{P(t)}^{1/t} = 
& \lim_{k\to\infty} \snorm{P(kT)}^{1/kT}
\\&
\leq (1 - \alpha^{T})^{1/T} 
\leq 1 - \alpha^{T} / T
\enspace.
\end{split}
\end{equation*}
Because all hypotheses we consider are time-invariant, it is sufficient to show
$S_j(T) = [n]$.

For Theorem~\ref{thm:coord:const}, we choose $T=\ind(H)$:
We show that $S_j(\ind(H)) = [n]$.
This is done by reducing the problem to one with a constant matrix.
So let~$A$ be any stochastic matrix whose digraph~$\graph{A}$ is equal to~$H$.
If~$A^t$ has a positive column, then so does~$P(t)$ because~$H$ is a
sub-digraph of every communication digraph.
This shows the claim since $\ind(\graph{A}) = \ind(H)$.

%
%

\section{Clusterings}\label{sec:clusterings}
We pair the idea of the distributed safety net in form of an aperiodic core
with the notion of {\em
clusters}, which have a leader that is the sole agent of the cluster to regard
values of agents other than the cluster's.
We will prove that it is not necessary for every agent to be contained in an
aperiodic component, but only for the cluster leaders.

A natural example of these clusterings occurs in the reduction of $B$-bounded
settings with self-confidence to synchronous ones (see
Fig.~\ref{fig:delay:reduction}), for which
$\ind(H)=B-1$.
If we do not assume self-confidence in $B$-bounded settings, then asymptotic
consensus is not necessarily reached, even if the averaging matrices are
constant and ergodic.
By proving results on cluster-aperiodic cores in synchronous settings, we are
hence also proving results on $B$-bounded settings with self-confidence.

\begin{figure}
\centering
\begin{tikzpicture}[>=latex',scale=0.7, every node/.style={transform shape}]
\node[draw,circle] (n1) at (0,3) {$1$};
\node[draw,circle] (n1p) at (-2,3) {$1'$};
\node[draw,circle] (n2) at (3,3) {$2$};
\node[draw,circle] (n3) at (3,0) {$3$};
\node[draw,circle] (n3p) at (5,0) {$3'$};
\node[draw,circle] (n3pp) at (7,0) {$3''$};
\node[draw,circle] (n4) at (0,0) {$4$};
\draw[->] (n4) -- (n1);
\draw[->] (n1) -- (n3);
\draw[->] (n3) .. controls +(-1.5,-0.4) .. (n4);
\draw[->] (n4) .. controls +(1.5,0.4) .. (n3);
\draw[->] (n2) .. controls +(-0.5,1.0) and +(0.5,1.0) .. (n2);
\draw[->] (n1p) -- (n1);
\draw[->] (n3p) -- (n3);
\draw[->] (n3pp) -- (n3p);
\node[dashed,cloud, cloud puffs=24, draw,minimum width=4cm, minimum
height=2cm] at (-1,3) {};
\node[dashed,cloud, cloud puffs=24, draw,minimum width=2cm, minimum
height=2cm] at (3,3) {};
\node[dashed,cloud, cloud puffs=24, draw,minimum width=2cm, minimum
height=2cm] at (0,0) {};
\node[dashed,cloud, cloud puffs=24, draw,minimum width=6cm, minimum
height=2cm] at (5,0) {};
\end{tikzpicture}
\caption{$\mathcal{C}$-aperiodic digraph with leaders $1,2,3,4$}
\label{fig:aperiodic:digraph}
\end{figure}
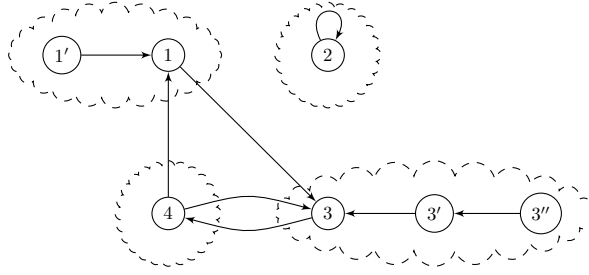

A digraph is a {\em cluster\/} with leader~$l$ if it is
$l$-coordinated.
A {\em clustering\/}~$\mathcal{C}$ is a collection of node-disjoint
clusters
$C_1,C_2,\dots,C_m$ together with respective leaders $l_1,l_2,\dots,l_m$.
A digraph is {\em
$\mathcal{C}$-aperiodic\/} if
every cluster~$C_j$ is a sub-digraph, every node is contained in some cluster,
and
it is $l$-aperiodic for every leader~$l_j$ of~$\mathcal{C}$.
Fig.~\ref{fig:aperiodic:digraph} shows an example of a
$\mathcal{C}$-aperiodic digraph.

A digraph {\em respects\/} a clustering~$\mathcal{C}$ if
the only edges leaving a cluster are the leader's.
Given a digraph that respects clustering~$\mathcal{C}$, the corresponding {\em
cluster digraph\/} is the digraph when collapsing all clusters of~$\mathcal{C}$
to single node.

\section{Dynamic Coordinated Communication Digraphs}\label{sec:dyn:coord}

We now prove that asymptotic consensus is also reached if there is no
coordinated core, but that coordination at every time step suffices.

\begin{thm}\label{thm:rooted}
A synchronous setting with averaging matrices $A(1),A(2),\dots$ 
with a $\mathcal{C}$-aperiodic 
spanning core~$H$ and minimal
confidence~$\alpha$ reaches asymptotic consensus if
every communication digraph respects~$\mathcal{C}$ and is coordinated.
Moreover, the rate of convergence is at most
\begin{equation*}
1 - \alpha^{(n-1)^2 (\ind(H)+1)}/(n-1)^2 (\ind(H)+1)
\end{equation*}
where~$n$ is the number of clusters in~$\mathcal{C}$.
\end{thm}

\begin{cor}\label{cor:rooted:delays}
A $B$-bounded setting with averaging matrices $A(1),A(2),\dots$ with
self-confidence and minimal
confidence~$\alpha$ reaches asymptotic consensus if
every communication digraph is coordinated.
Moreover, the rate of convergence is at most
$1 - \alpha^{(n-1)^2 B}/(n-1)^2 B$.
\end{cor}

Corollary~\ref{cor:rooted:delays}, without the explicit bound on the rate of
convergence is included in Theorem~\ref{thm:cao:et:al}.

We prove the theorem in the rest of the subsection.
Recall that~$n$ is the number of nodes, not the number of clusters.
We also note that $T(C_j)\leq T(H)$ for every cluster~$C_j$ in~$\mathcal{C}$.

The sets~$S_{j}(t)$ 
satisfy a weak form of monotonicity if the sequence of communication
graphs have an aperiodic core.
If there are self-loops in all communication digraphs, then clearly
$S_j(t)\subseteq S_j(t+1)$, which is a special case of the following lemma.

\begin{lem}\label{lem:sj:monotone}
If~$H$ is a spanning $\mathcal{C}$-aperiodic core and all communication digraphs
respect~$\mathcal{C}$,
then 
$S_j(t_1) \subseteq S_j(t_2)$ whenever $t_2 - t_1\geq \ind(H)$ and~$j$ is a
leader of~$\mathcal{C}$.
\end{lem}
\begin{proof}
Let $i\in S_j(t_1)$.
Since all communication digraphs respect the clustering, $i$'s leader~$l_i$
appears in some earlier set: $l_i\in S_j(t_1')$ with $t_1'\leq t_1$.

Because~$H$ is $l_i$-aperiodic and $t_2-t_1'\geq\ind(H)$,
there exists
a walk of length~$t_2-t_1'$ from~$i$ to~$l_i$ in~$H$
by the definition of $\ind(H)$.
The fact that~${H}$ is a sub-digraph of all~$G\big(A(\tau)\big)$
shows that $P_{i,l_i}(t_2,t_1')$ is positive by Lemma~\ref{lem:graph}.

Hence
\begin{equation*}
P_{i,j}(t_2) = 
\sum_k P_{i,k}(t_2,t_1') \cdot P_{k,j}(t_1')
\geq
P_{i,l_i}(t_2,t_1') \cdot P_{l_i,j}(t_1')
\end{equation*}
is positive, which shows $i\in S_j(t_2)$.
\end{proof}

The following lemmas are used to lower bound the steps need until $S_j(t)=[n]$.

\begin{lem}\label{lem:roots:and:sj}
If~$H$ is a spanning $\mathcal{C}$-aperiodic core, all communication graphs
respect~$\mathcal{C}$, $j$ is a leader of~$\mathcal{C}$, $t\geq \ind(H)$, and $\graph{A(t+1)}$ is
$j$-coordinated, then 
either $S_j(t) = [n]$ or 
$S_j(t+1)\setminus S_j(t) \neq \emptyset$.
\end{lem}
\begin{proof}
The hypothesis that $t\geq \ind(H)$ guarantees that $j\in S_j(t)$ by
Lemma~\ref{lem:sj:monotone}.
Every node has a path to~$j$,  and hence
to $S_j(t)$,  in $\graph{A(t+1)}$.
Now, if $S_j(t)\neq[n]$, there is some $i\in[n]\setminus S_j(t)$ that has an
outgoing neighbor~$k_0$ in~$S_j(t)$, i.e.,
$A_{i,k_0}(t+1)>0$.
The condition $k_0\in S_j(t)$ means $P_{k_0,j}(t)>0$ and hence
\begin{equation*}
\begin{split}
P_{i,j}(t+1) & = \sum_{k} A_{i,k}(t+1) \cdot P_{k,j}(t-1)
\\ &
\geq
A_{i,k_0}(t+1) \cdot P_{k_0,j}(t)>0
\enspace,
\end{split}
\end{equation*}
which shows $i\in S_j(t+1)$.
\end{proof}

\begin{lem}\label{lem:cluster:follows:leader}
Let~$H$ be a spanning $\mathcal{C}$-aperiodic core, all communication graphs
respect~$\mathcal{C}$ and $j$ be a leader of~$\mathcal{C}$.
If~$l$ is any leader of some cluster~$C$ of~$\mathcal{C}$ and $l\in
S_j(t)$, then $C\subseteq S_j(t + \ind(H))$.
\end{lem}
\begin{proof}
Because $C$ is $l$-aperiodic and $l$-coordinated, we have $C\subseteq
S_{l}(\tau)$ for all $\tau \geq \ind(C)$.
Because $\ind(H)\geq\ind(C)$, the lemma follows with an application of
Lemma~\ref{lem:graph}.
\end{proof}

Set $t_m = m\cdot (\ind(H) + 1)$.
For $m\geq1$, let~$j_m$ be a leader of the digraph~$G(A(t_m))$ and also
of~$\mathcal{C}$.
Lemma~\ref{lem:sj:monotone} specialized to $s=t_{m-1}$ and $t = t_m - 1 =
t_{m-1} + \ind(H)$ gives
$S_j(t_m - 1) \supseteq S_j(t_{m-1})$
for all leaders~$j$ and all $m\geq1$.
Lemma~\ref{lem:roots:and:sj} applied to $t = t_m$ and $j=j_m$ gives:
$S_{j_m}(t_m) \supsetneq S_{j_m}(t_{m-1})$
if $S_{j_m}(t_m - 1) \neq [n]$.

If $m= (n-1)^2 = (n-2)n + 1$, then some $j_0\in[n]$ appears at least $n-1$
times in the sequence of leaders $j_1,j_2,\dots,j_m$.
By the above and Lemma~\ref{lem:cluster:follows:leader}, it is hence
$S_{j_0}(t_m) = [n]$,
which shows the theorem. 

\section{Dynamic Communication Digraphs with Fixed Leader}\label{sec:dyn:fixed}

In this subsection, we assume a {\em fixed\/} leader in every communication digraph
and are able to show a tighter bound on the rate of convergence.
The case of strongly connected communication digraphs is a special case.


\begin{thm}\label{thm:fixed:root}
A synchronous setting with averaging matrices $A(1),A(2),\dots$ 
with a $\mathcal{C}$-aperiodic 
spanning core~$H$ and minimal
confidence~$\alpha$ reaches asymptotic consensus if
every communication digraph respects~$\mathcal{C}$ and
there is an agent~$j_0$ such that every communication digraph is
$j_0$-coordinated.
Moreover, the rate of convergence is at most
\begin{equation}
1 - \alpha^{(n-1)(\ind(H)+1)}/(n-1)(\ind(H)+1)
\end{equation}
where~$n$ is the number of clusters in~$\mathcal{C}$.
\end{thm}

\begin{cor}\label{cor:rooted:fixed:delays}
A $B$-bounded setting with averaging matrices $A(1),A(2),\dots$ with
self-confidence and minimal
confidence~$\alpha$ reaches asymptotic consensus if
there is an agent~$j_0$ such that
every communication digraph is $j_0$-coordinated.
Moreover, the rate of convergence is at most
$1 - \alpha^{(n-1) B}/(n-1) B$.
\end{cor}

Corollary~\ref{cor:rooted:fixed:delays}, without the explicit bound on the rate
of convergence is included in Theorem~\ref{thm:cao:et:al}.

%
%

We use the notation of the previous subsection.
The theorem follows similarly by noticing that, in this case, $j_m=j_0$ for all
$m\geq
1$ and hence~$j_0$ appears $n-1$ times in the sequence of leaders
$j_1,j_2,\dots,j_{n-1}$.

\section{Completely Reducible Communication Digraphs}\label{sec:red}

We now show that one can replace the assumption of coordination by the
assumption of completely reducibility at every time step and eventual weak
connectivity.

\begin{thm}\label{thm:aperiodic}
A synchronous setting with averaging matrices $A(1),A(2),\dots$ 
with a $\mathcal{C}$-aperiodic 
spanning core~$H$ and minimal
confidence~$\alpha$ reaches asymptotic consensus if
every communication digraph respects~$\mathcal{C}$,
all cluster communication digraphs are completely reducible, and
the digraph~$G_\infty$ formed by all edges that appear in infinitely many
cluster communication digraphs is weakly connected.
\end{thm}

\begin{cor}\label{cor:reducible:delays}
A $B$-bounded setting with averaging matrices $A(1),A(2),\dots$ with
self-confidence and minimal
confidence~$\alpha$ reaches asymptotic consensus if
every communication digraph is completely reducible and the digraph~$G_\infty$
of edges that appear in infinitely many communication digraphs is weakly
connected.
\end{cor}

Corollary~\ref{cor:reducible:delays}
for synchronous settings is  Theorem~\ref{thm:touri:nedic}.

We prove this theorem in the rest of this subsection.
%
%
%
We do not use the exact same proof strategy as in the previous subsection:
We show the existence of a~$T$ such that 
\begin{equation*}
\snorm{ P(T)} \leq  1 - \alpha^{n(\ind(H)+1)}
\enspace.
\end{equation*}
This suffices to show the theorem because the conditions in the theorem are
time-invariant and
repeated application thus shows that $\snorm{ P(t)} \to 0$.
Even though we cannot bound~$T$ with the hypotheses of the theorem, we {\em
can\/} bound the semi-norm uniformly, which is critical for the proof to work.
Lemma~\ref{lem:norm:limit} then concludes the proof.

We first show that~$G_\infty$ is completely reducible.
For that, we show the following basic lemma.

\begin{lem}
\label{lem:compl:red:stable}
Every union of completely reducible digraphs
is completely reducible.
\end{lem}
\begin{proof}
Let~$\mathcal{G}$ be a set of completely reducible digraphs and let $H =
\bigcup \mathcal{G}$ be their union.
Let~$i$ and~$j$ be two nodes in~$H$ and suppose that there exists a path~$P$
from~$i$ to~$j$ in the union digraph~$H$.
We will show that there then exists a path from~$j$ to~$i$ in~$H$.
This is trivial if~$i=j$ so suppose the contrary, i.e., that~$P$ is nonempty.

Let $i_0,i_1,\dots,i_n$ be $P$'s sequence of nodes.
For every $1\leq k\leq n$, the edge~$e_k$ is in some digraph $G\in\mathcal{G}$.
Now, because~$G$ is completely reducible, there exists a path~$P_k$ in~$G$
from~$e_k$ to~$e_{k-1}$.
But then the composite walk $P_n\cdot P_{n-1}\cdots P_1$ is a walk in~$H$
from~$j$ to~$i$.
\end{proof}


Hence~$G_\infty$ is completely reducible because
Lemma~\ref{lem:compl:red:stable} shows that
\begin{equation*}
G_\infty = \lim_{T\to\infty} \bigcup_{t\geq T} G(A(t))
\end{equation*}
is a decreasing limit of a sequence of completely reducible digraphs.
Because all digraphs are finite, this sequence is eventually constant.
Hence its limit~$G_\infty$ is equal to one of the sequence's elements and hence
completely reducible.

The next lemma captures the essence of the complete reducibility assumption:
If $S_j(t)$ does not change, then $\mu_j(t)$ does not decrease.
Together with the weak monotonicity of Lemma~\ref{lem:sj:monotone} and eventual
connectivity, we are able to show the theorem.

\begin{lem}\label{lem:equal}
Under the hypotheses of Theorem~\ref{thm:aperiodic}, if $j$ is a leader
of~$\mathcal{C}$ and
$S_{j}(t)=S_{j}(t+1)$, then $\mu_{j}(t+1)\geq \mu_{j}(t)$.
\end{lem}
\begin{proof}
Let $P_{i,{j}}(t+1)$ be positive, i.e., $i\in S_{j}(t+1)=S_{j}(t)$.
By definition of~$S_{j}(t)$, we have
\begin{equation}\label{eq:proof:equal:1}
P_{i,j}(t+1) = \sum_{k\in S_{j}(t)} A_{i,k}(t+1) \cdot P_{k,j}(t)
\enspace.
\end{equation}
Because $S_{j}(t)=S_{j}(t+1)$, we derive that $A_{i,k}(t+1)$ is zero
whenever
$i\not\in
S_{j}(t)$ and $k\in S_{j}(t)$.
Because every node of a cluster is leader-coordinated, every the nodes of a
cluster are either all in~$S_j(t)$ or all outside of~$S_j(t)$.
Hence, because the cluster digraph~$A(t+1)$ is completely irreducible, we also
have that
$A_{i,k}(t+1)$ is zero whenever $i\in
S_{j}(t)$ and $k\not\in S_{j}(t)$.

By assumption, we have $i\in S_{j}(t)$, and hence by the above and by
stochasticity of~$A(t+1)$:
\begin{equation}\label{eq:proof:equal:2}
1 = \sum_{k} A_{i,k}(t+1) = \sum_{k\in S_{j}(t)} A_{i,k}(t+1)
\end{equation}
Because $P_{k,j}(t)\geq\mu_{j}(t)$ for all $k\in S_{j}(t)$, combination
of Equations~\eqref{eq:proof:equal:1} and~\eqref{eq:proof:equal:2} yields
$P_{i,j}(t+1)\geq\mu_{j}(t)$.
\end{proof}


Choose any leader $j_0$ of~$\mathcal{C}$.
For every~$i\in[n]$, let~$t_i$ be the least nonnegative integer such that
$C_i\subseteq
S_{j_0}(t_i)$.
All~$t_i$ are
well-defined as~$G_\infty$ is strongly connected.
By permuting indices, we can assume without loss of generality that $t_1\leq t_2 \leq \cdots \leq
t_n$.
Because~$P(0)$ is the identity matrix, we have $ S_{j_0}(0) =  \{j_0\}$ and hence $t_1=0$.

We inductively show
\begin{equation}\label{eq:proof:lem:final:mu}
\mu_{j_0}(t_m) \geq \alpha^{(m - 1) (\ind(H)+1)}
\end{equation}
for
all $1\leq m\leq n$.
This is true for~$m=1$.
To prove the inductive step, we distinguish two cases:
(A) $t_{m}-t_{m-1} < \ind(H)$ and (B) $t_{m}-t_{m-1} \geq \ind(H)$.

In case~(A),
we have
\begin{align*}
\mu_{j_0}(t_m)& \geq \alpha^{t_m-t_{m-1}} \cdot \mu_{j_0}(t_{m-1}) \geq
\alpha^{(m-1)(\ind(H)+1)}
\end{align*}
by the induction hypothesis.

In case~(B), we have $S_{j_0}(t)=S_{j_0}(t_{m-1})$ for all~$t$ with
$t_{m-1}+\ind(H)\leq t\leq t_{m}-1$
by Lemma~\ref{lem:sj:monotone} and the definition of~$t_m$.
Repeated application of Lemma~\ref{lem:equal} hence yields
$\mu_{j_0} (t_m-1)\geq \mu_{j_0}\big(t_{m-1}+\ind(H)\big)$.
We thus have
\begin{equation*}
\begin{split}
\mu_{j_0}(t_m) 
& \geq 
\alpha \cdot \mu_{j_0}(t_m-1) 
\geq 
\alpha \cdot \mu_{j_0}(t_{m-1}+\ind(H))
\\ & \geq
\alpha^{\ind(H)+1} \cdot \mu_{j_0}(t_{m-1})
\geq
\alpha^{(m-1)(\ind(H)+1)}
\end{split}
\end{equation*}
by the induction hypothesis.

In particular, we have shown Equation~\eqref{eq:proof:lem:final:mu} for
$m=n$.
Now set $T=t_n+\ind(H)$.
By Lemmas~\ref{lem:sj:monotone} and~\ref{lem:cluster:follows:leader},
$S_{j_0}(T)=[n]$ for all
and
$\mu_{j_0}(T) \geq \alpha^{n(\ind(H)+1)}$.
This concludes the proof of the theorem.

\section{Conclusion}\label{sec:conclusion}
The paper introduced the novel notion of aperiodic cores and showed that the
prevalent hypothesis of self-confidence can be replaced by the hypothesis of
the existence of an aperiodic core in a large variety of convergence results
for asymptotic consensus in dynamic settings.
In particular, we discussed and explored the case of non-synchronous
environments, for which we gave an explicit example of a $2$-bounded system
with $2$ agents that could not be handled by existing convergence theorems.
We also highlighted the need to be careful in these matters by showing that a
small variant of the example does not reach asymptotic consensus (and does not
even converge).
In a linear algebraic view, our results are strict generalizations of the
Perron-Frobenius theorem, which was not the case for most results on asymptotic
consensus in the literature, as they require self-confidence.

\end{document}